\documentclass[twocolumn]{svjour3}
\usepackage{array,amsmath,amstext,amssymb,amsbsy,graphicx}
\usepackage{stmaryrd}
\usepackage{pdfsync}
 \usepackage{color,bm}
\usepackage{txfonts}
\usepackage{todonotes}

\newcommand{\nablas}{\nabla_\Gamma}
\newcommand{\nablash}{\nabla_{\Gamma_h}}
\newcommand{\mean}[1]{\langle {#1} \rangle}
\newcommand{\jump}[1]{[#1]}

\newcommand{\IR}{\mathbb{R}}

\newcommand{\bfJ}{{\boldsymbol J}}

\newcommand{\bfU}{{\boldsymbol U}}

\newcommand{\bfn}{\boldsymbol n}

\newcommand{\bfkappa}{\boldsymbol\kappa}

\newcommand{\bfF}{\boldsymbol F}
\newcommand{\bftheta}{\boldsymbol\theta}
\newcommand{\bfv}{\boldsymbol v}
\newcommand{\bfxi}{\boldsymbol\xi}

\newcommand{\bfgamma}{\boldsymbol\gamma}

\newcommand{\bff}{\boldsymbol f}
\newcommand{\bfp}{\boldsymbol p}
\newcommand{\bfu}{\boldsymbol u}
\newcommand{\bfeps}{\boldsymbol\varepsilon}
\newcommand{\bfe}{\boldsymbol e}
\newcommand{\bfsig}{\boldsymbol\sigma}
\newcommand{\bfnu}{\boldsymbol\nu}
\newcommand{\bftau}{\boldsymbol\tau}

\newcommand{\bfw}{\boldsymbol w}
\newcommand{\bfI}{\boldsymbol I}
\newcommand{\bfP}{\boldsymbol P}
\newcommand{\bfx}{\boldsymbol x}

\newcommand{\bfV}{\boldsymbol V}

\newcommand{\mcK}{\mathcal{K}}
\newcommand{\mcE}{\mathcal{E}}

\newcommand{\mcX}{\mathcal{X}}
\newcommand{\mcI}{\mathcal{I}}

\newcommand{\ttt}{\tilde{t}^2}
\numberwithin{equation}{section}

\newtheorem{rem}{Remark}[section]

\begin{document}
\title{Continuous/Discontinuous Finite Element Modelling of Kirch\-hoff Plate 
Structures in $\mathbb{R}^3$ Using Tangential Differential Calculus.}
\author{Peter Hansbo \and Mats G. Larson}
\institute{P. Hansbo\newline Department of Mechanical Engineering, J\"onk\"oping University,
SE-55111 J\"onk\"oping, Sweden \newline\newline M. G. Larson \newline Department of Mathematics and Mathematical Statistics,
Ume{\aa} University, SE--901 87 Ume{\aa}, Sweden}
%
\date{\today}

\maketitle

\begin{abstract}
We employ surface differential calculus to derive models for Kirchhoff plates including in--plane membrane deformations. We also extend our formulation to structures of plates. For solving the resulting set of partial differential equations, we employ a finite element method based on elements that are continuous for the displacements and discontinuous for the rotations, using $C^0$--elements for the discretisation of the plate as well as for the membrane deformations. Key to the formulation of the method is a convenient definition of jumps and averages of forms that are $d$-linear in terms of the element edge normals.
\end{abstract}

\section{Introduction}

The Kirchhoff plate model is a fourth order partial differential equation which requires $C^1$--continuous elements for constructing
conforming finite element methods. To avoid this requirement, nonconforming finite elements can be used; one classical example being the Morley triangle \cite{Mo68} which has displacement degrees of freedom in the corner nodes and rotation degrees of freedom at the midpoint of the edges. If we want to solve also for the membrane displacements, 
it is more straightforward to be able to use only displacement degrees of freedom for both the normal (plate) and tangential (membrane) displacements.
To reach this goal, one can instead use the discontinuous Galerkin (dG) method \cite{HaLa02}, more efficiently implemented as a $C^0$--continuous Galerkin method 
allowing for discontinuous approximation of derivatives, referred to as the continuous/discontinuous Galerkin, or c/dG, method, 
first suggested by Engel et al. \cite{engel02}, and further developed for
plate models by Hansbo et al. \cite{HaHeLa10,HaHeLa11,HaLa03,HaLa11} and by Wells and Dung \cite{WeDu07}. See also Larsson and Larson \cite{LaLa17} for error estimates in the case of the biharmonic problem on a surface.
To obtain a continuous model, we combine the plate equation for the normal displacements with the tangential differential equation for the membrane from Hansbo and Larson \cite{HaLa14a} to obtain a structure with both bending resistance and membrane action. This model is then discretised using continuous finite elements for the membrane and c/dG for the plate, using the same order polynomial in both cases. 

The standard engineering approach to constructing plate elements arbitrarily oriented in ${\mathbb{R}}^3$ is to
use rotation matrices to transform the displacements from a planar element to the actual, common, coordinates, thus transforming the stiffness matrices. In this paper we instead extend the c/dG
method to the case of arbitrarily oriented plates, allowing for membrane deformations, directly using Cartesian coordinates in ${\mathbb{R}}^3$. We argue that this makes it simpler to implement discrete schemes in general, and in particular the discontinuous Galerkin terms on the element borders. It also gives an analytical model directly expressed in equilibrium equations in physical coordinates.

A particular feature of our method is the handling of the trace terms in the c/dG method. In the recent paper on dG for elliptic problems on smooth surfaces by Dedner, Madhavan, and Stinner \cite{DeMaSt13} the definition of the normal to the element faces (tangential to the surface), the \emph{conormal},\/ was discussed and different variants tested numerically. In our case, where the surface is piecewise smooth (planar), the definition of the conormal at plate junctures is crucial to the equilibrium. It turns out the 
proper way to define the jumps and averages of trace quantities that are $d$-linear in the conormal 
is to compute the trace on the left and right side with the respective unit conormals and adjust the 
sign on one of the sides with $(-1)^d$. This leads to a generalization of the standard jump and averages 
in the flat case where a fixed conormal is used for both the left and right side in the definition of the jump. Furthermore, the standard formula, where the jump in a product of two functions is represented as the 
sum of the two products of the averages and jumps of the two functions, also generalizes to this situation.
With these tools at hand we may directly use standard discontinuous Galerkin techniques to 
derive a finite element method for a plate structure. The resulting method takes the same form as a 
standard c/dG method for a plate. The only difference is the proper definition of jumps and averages.
See also \cite{JoLaLa17}, where a similar approach was used for the Laplace-Beltrami operator on a surface with sharp edges.

The outline of the paper is as follows: In Section 2 we derive a variational formulation for a plate 
with arbitrary orientation in $\IR^3$, in Section 3 we define the relevant traces, including forces and 
moments, define the averages and jumps of $d$-linear forms, and formulate the interface conditions 
for a plate structure, in Section 4 we formulate the finite element method, in Section 5 we present numerical examples, and finally we conclude with some remarks in Section 6.

\section{Single Plate}

\subsection{Tangential Differential Calculus}
Let $\Gamma$ be a piecewise planar two-dimensional surface imbedded in ${\mathbb{R}}^3$, with piecewise constant unit normal $\bfn$ and boundary $\partial\Gamma$, split into a Neumann part $\partial\Gamma_\text{N}$ where forces and moments are known, and a Dirichlet part $\partial\Gamma_\text{D}$ where rotations and displacements are known. 
For ease of presentation we shall assume that $\partial\Gamma_\text{N}=\emptyset$ and that we have zero displacements and rotations on the boundary.
The case of  $\partial\Gamma_\text{N}\neq\emptyset$ is straightforward to implement and will be used in the numerical examples.  Mixed boundary conditions are handled equally straightforward.

If we denote the (piecewise) signed distance function relative to $\Gamma$ by $\zeta(\bfx)$, for $\bfx\in \mathbb{R}^3$, fulfilling $\nabla\zeta = \bfn$, we can define the domain occupied by the shell by
\begin{equation}
\Omega_t = \{\bfx\in \mathbb{R}^3: \vert \zeta(\bfx) \vert < t/2\}
\end{equation}
where $t$ is the thickness of the shell, which for simplicity will be assumed constant. The closest point projection $\bfp:\Omega_t \rightarrow \Gamma$
is given by
\begin{equation}
\bfp(\bfx) = \bfx -\zeta(\bfx)\bfn(\bfx) 
\end{equation}
the Jacobian matrix of which is 
\begin{equation}
\nabla \bfp = \bfI -\zeta\nabla\otimes\bfn -\bfn\otimes\bfn
\end{equation}
where $\bfI$ is the identity and $\otimes$ denotes exterior product.
The corresponding linear projector
$\bfP_\Gamma = \bfP_\Gamma(\bfx)$, onto the tangent plane of $\Gamma$ 
at $\bfx\in\Gamma$, is given by
\begin{equation}
\bfP_\Gamma := \bfI -\bfn\otimes\bfn 
\end{equation}
and we can then define the surface gradient
$\nablas$ as
\begin{equation}
\nablas := \bfP_\Gamma \nabla 
\end{equation}
The surface gradient thus has three
components, which we shall denote by
\begin{equation}
\nablas =: \left[\begin{array}{>{\displaystyle}c}
\frac{\partial}{\partial x_\Gamma} \\[3mm]
\frac{\partial}{\partial y_\Gamma} \\[3mm]
\frac{\partial}{\partial z_\Gamma} 
\end{array}\right] 
\end{equation}
For a vector valued function $\bfv(\bfx)$, we define the
tangential Jacobian matrix as
\begin{equation}
\bfv\otimes \nablas :=\left[\begin{array}{>{\displaystyle}c>{\displaystyle}c>{\displaystyle}c}
\frac{\partial v_1}{\partial x_\Gamma} &\frac{\partial v_1}{\partial y_\Gamma} & \frac{\partial v_1}{\partial z_\Gamma} \\[3mm]
\frac{\partial v_2}{\partial x_\Gamma} &\frac{\partial v_2}{\partial y_\Gamma} & \frac{\partial v_2}{\partial z_\Gamma} \\[3mm]
\frac{\partial v_3}{\partial x_\Gamma} &\frac{\partial v_3}{\partial y_\Gamma} & \frac{\partial v_3}{\partial z_\Gamma}
\end{array}\right] 
\end{equation}
and the surface divergence $\nabla_{\Gamma}\cdot\bfv := \text{tr}\,\bfv\otimes \nablas$.

\subsection{Displacement and Strain}
Upon loading, each point $\bfx \in \Omega_t$, in the plate
undergoes a displacement
\begin{equation}
\bfu (\bfx)  =  \bfu_0 (\bfp(\bfx))-\zeta(\bfx)\bfw(\bfp(\bfx))
\end{equation}
where $\bfu_0$ and $\bfw$ are vector fields defined on $\Gamma$, $\bfu_0$ 
arbitrary and $\bfw$ a tangential vector, $\bfw\cdot\bfn=0$ on $\Gamma$, or 
$\bfw = \bfP_\Gamma \bftheta$ with $\bftheta$ arbitrary.
Thus, neglecting in-plane extensions for the moment, we can write
\begin{equation}
\bfu = u_n \bfn-\zeta\bfP_\Gamma\bftheta
\end{equation}
in $\Omega_t$. Here $u_n = \bfu \cdot \bfn$. 

We introduce the strain tensor 
$\bfeps$ as
\begin{equation}
\bfeps(\bftheta) :=\frac{1}{2}\left(\bftheta\otimes \nabla + (\bftheta\otimes \nabla)^{\rm T}\right) 
\end{equation}
and define the symmetric part of the tangential Jacobian as
\begin{equation}
\bfe_{\Gamma}(\bftheta) := \frac12\left(\bftheta\otimes \nablas + (\bftheta\otimes \nablas)^{\rm T}\right) 
\end{equation}
The in-plane strain tensor $\bfeps_\Gamma$ is implemented using the 
following identity 
\begin{align}
\bfeps_\Gamma(\bftheta) 
&=  \bfP_\Gamma\bfe(\bftheta)\bfP_\Gamma 
\\ \label{eq:identity-a}
&=  \bfe_{\Gamma}(\bftheta) - (\bfe_{\Gamma}(\bftheta)\cdot\bfn)\otimes\bfn  
- \bfn\otimes (\bfe_{\Gamma}(\bftheta)\cdot\bfn) 
\end{align}
If we write
\begin{equation}
\bftheta = \bfP_\Gamma\bftheta + (\bftheta\cdot\bfn)\bfn
\end{equation}
then
\begin{equation}
\bfeps_{\Gamma}(\bftheta) = \bfeps_{\Gamma}( \bfP_\Gamma\bftheta) +(\bftheta\cdot\bfn)\bfkappa
\end{equation}
where
\begin{equation}
\bfkappa := \nabla\otimes\bfn
\end{equation}
is the curvature tensor,
cf. \cite{DeZo95,DeZo96}. For planar $\Gamma$, $\bfn$ is constant, and this 
simplifies to 
\begin{equation}\label{planar}
\bfeps_{\Gamma}(\bftheta) = \bfeps_{\Gamma}( \bfP_\Gamma\bftheta)  
\end{equation}
The total in-plane strain tensor is thus given by
\begin{equation}
\bfeps_{\Gamma}(\bfu) = \bfeps_{\Gamma}(u_n \bfn)-\zeta \bfeps_{\Gamma}(\bfP_\Gamma\bftheta) 
\end{equation}

%
%
In \cite{DeZo95,DeZo96} it is also shown that the mid-plane rotation in the absence of shear deformation is given by $2\bfe_{\Gamma}(u_n \bfn)\cdot\bfn$,
and for shear deformable inextensible shells we thus have the shear deformation vector
\begin{equation}
\bfgamma = \frac12 \left(2\bfe_{\Gamma}(u_n\bfn)\cdot\bfn - \bfP_\Gamma\bftheta\right)
\end{equation}
It is is also easy to verify that
\begin{equation}
\bfe_{\Gamma}(u_n\bfn) = u_n\,\bfe_\Gamma(\bfn)+\frac12\left(\bfn\otimes\nablas u_n+ (\bfn\otimes\nablas u_n)^{\rm T}\right)
\end{equation}
so that, since $\bfn\cdot\nablas u_n  = 0$,
\begin{equation}
2\bfe_{\Gamma}(u_n \bfn)\cdot\bfn = \nablas u_n +2 \bfe_\Gamma(\bfn) \cdot\bfn \, u_n =  \nablas u_n 
\end{equation}
since $\bfn$ is constant; thus
\begin{equation}
\bfgamma = \frac12 \left(\nablas u_n - \bfP_\Gamma\bftheta\right)
\end{equation}
In the tangential setting, the Kirchhoff assumption of zero shear deformations can therefore 
be written
\begin{equation}
\bfu := u_n\bfn - \zeta\nablas u_n 
\end{equation}
Furthermore, we find that
\begin{align}
\bfeps_\Gamma(u_n \bfn)  &= \bfP_\Gamma \bfe_\Gamma(u_n \bfn ) \bfP_\Gamma
\\
&= \frac{1}{2} \bfP_\Gamma  \Big( ( \nablas u_n) \otimes \bfn + \bfn \otimes (\nablas u_n) \Big) \bfP_\Gamma
\\
&= 0
\end{align}
and for inextensible plates we get
\begin{equation}
\bfeps_\Gamma(\bfu) = - \zeta \bfeps_\Gamma(\nablas u_n )
\end{equation}
and in this case we thus only obtain contributions to the strain energy 
from the displacement field
\begin{equation}
\bfu  = -\zeta\nablas u_n 
\end{equation}

\subsection{Variational Formulations}

We shall assume isotropic stress--strain relations,
\begin{equation}\label{eq:Hooke}
\bfsig = 2\mu \bfeps + \lambda \text{tr}\,\bfeps \, \bfI 
\end{equation}
where $\bfsig$ is the  stress tensor, and plane stress conditions, 
for which the Lam\'e parameters $\lambda$ and $\mu$ are related to Young's modulus $E$ 
and Poisson's ratio $\nu$ via
\begin{equation}
\mu = \frac{E}{2(1+\nu)},\quad \lambda =\frac{E\nu}{1-\nu^2} 
\end{equation}
For the in-plane stress tensor we find, by projecting (\ref{eq:Hooke}) from left and 
right,
\begin{equation}\label{stresstrain}
\bfsig_\Gamma:=  2\mu \bfeps_\Gamma + {\lambda} \text{tr}\bfeps_\Gamma\, \bfP_\Gamma = 2\mu \bfeps_\Gamma + {\lambda} \nabla_\Gamma\cdot\bfu\, \bfP_\Gamma
\end{equation}
The potential energy of the plate is postulated as
\begin{align}
{\cal E}_\text{P}&:= \frac12\int_{-t/2}^{t/2}\int_\Gamma\bfsig_\Gamma(\zeta \nablas u) : \bfeps_\Gamma(\zeta \nablas u)\, d\Gamma d\zeta  
\\ \nonumber
&\quad -\int_{-t/2}^{t/2}\int_\Gamma\bff\cdot\bfu\, d\Gamma d\zeta
\end{align}
where $\bfsig :\bfeps = \sum_{ij}\sigma_{ij}\varepsilon_{ij}$ for second order Cartesian tensors $\bfsig$ and $\bfeps$.
Integrating in $\zeta$, we obtain
\begin{align}
{\cal E}_\text{P}  :=  & {} \frac{t^3}{24}\int_\Gamma{\bfsig}_\Gamma(\nablas u) : {\bfeps}_\Gamma(\nablas u)\, d\Gamma   -t\int_\Gamma\bff\cdot\bfn \, u\, d\Gamma 
\end{align}
Under the assumption of clamped boundary conditions, the corresponding variational 
problem is to find $u_n \in H_0^2(\Gamma) = \{\text{$v \in H^2(\Gamma)$ : 
$v = \bfnu \cdot \nablas v = 0$ on $\partial \Gamma$}\}$ such
that
\begin{equation}
\frac{t^3}{12}\int_\Gamma{\bfsig}_\Gamma(\nablas u) : {\bfeps}_\Gamma(\nablas v)\, d\Gamma   
 =t\int_\Gamma \bff\cdot\bfn v\, d\Gamma 
\end{equation}
for all $v\in H^2_0(\Gamma)$.

Introducing also membrane deformations, the total potential energy ${\cal E}_\text{tot}$ 
of the plate must take into account both the bending energy ${\cal E}_\text{P}$ and the membrane energy ${\cal E}_\text{M}$, so that 
${\cal E}_\text{tot}={\cal E}_\text{P}+{\cal E}_\text{M}$, 
where
\begin{align}
{\cal E}_\text{M}  
&:=  t\int_\Gamma{\bfsig}_\Gamma( \bfP_\Gamma\bfu_0) : {\bfeps}_\Gamma( \bfP_\Gamma\bfu_0)\, d\Gamma   
\\ \nonumber
&\quad -t \int_\Gamma\bff\cdot \bfP_\Gamma \bfu_0\, d\Gamma 
\end{align}

Since we wish to use a 3D Cartesian vector field we redefine $\bfu :=\bfu_0$ and $u_n:=\bfn\cdot\bfu$, make use of (\ref{planar}), and introduce the function space
\begin{equation}
V=\{\bfv: \;\bfP_\Gamma\bfv\in [H_0^1(\Gamma)]^3, \; v_n = \bfv\cdot \bfn \in H^2_0(\Gamma)\} .
\end{equation}
We are then led to the variational problem of finding $\bfu \in V$ such
that
\begin{align} \nonumber
&\frac{t^2}{12} \int_\Gamma{\bfsig}_\Gamma(\nablas u_n) : {\bfeps}_\Gamma(\nablas v_n)\, d\Gamma 
\\ \label{eq:variational-a}
&\quad + \int_\Gamma{\bfsig}_\Gamma(\bfu) : {\bfeps}_\Gamma(\bfv)\, d\Gamma   
 =   \int_\Gamma \bff\cdot \bfv\, d\Gamma 
\end{align}
for all $\bfv\in V$. Introducing the notation 
\begin{equation}
 \tilde{t} = \frac{t}{\sqrt{12}}
\end{equation}
we may write (\ref{eq:variational-a}) in the more compact form
\begin{align} \nonumber
&\ttt\int_\Gamma{\bfsig}_\Gamma(\nablas u_n) : {\bfeps}_\Gamma(\nablas v_n)\, d\Gamma 
\\
&\quad + \int_\Gamma{\bfsig}_\Gamma(\bfu) : {\bfeps}_\Gamma(\bfv)\, d\Gamma   
 =   \int_\Gamma \bff\cdot \bfv\, d\Gamma 
\end{align}

For implementation purposes we note that for $\bfn$ constant
\begin{equation}
\nablas u_n = (\bfu\otimes\nablas)\cdot\bfn
\end{equation}
and 
\begin{equation}
\begin{array}{>{\displaystyle}c}
\bfe_{\Gamma}\left((\bfu\otimes\nablas)\cdot\bfn\right) =\\[3mm]
\left[\begin{array}{>{\displaystyle}c>{\displaystyle}c>{\displaystyle}c}
\frac{\partial^2 \bfu}{\partial x_\Gamma^2}\cdot\bfn &\frac{\partial^2 \bfu}{\partial x_\Gamma\partial y_\Gamma}\cdot\bfn & \frac{\partial^2\bfu}{\partial x_\Gamma\partial z_\Gamma}\cdot\bfn \\[3mm]
\frac{\partial^2 \bfu}{\partial x_\Gamma\partial y_\Gamma}\cdot\bfn &\frac{\partial^2 \bfu}{\partial y_\Gamma^2}\cdot\bfn & \frac{\partial^2 \bfu}{\partial y_\Gamma\partial z_\Gamma}\cdot\bfn \\[3mm]
\frac{\partial^2\bfu}{\partial x_\Gamma\partial z_\Gamma}\cdot\bfn  &\frac{\partial^2 \bfu}{\partial y_\Gamma\partial z_\Gamma}\cdot\bfn  & \frac{\partial^2 \bfu}{\partial z_\Gamma^2}\cdot\bfn
\end{array}\right] \end{array}
\end{equation}

\subsection{Strong Form}

The corresponding strong form of the problem is to find $\bfu = u_n\bfn + \bfP_\Gamma\bfu$ such that
\begin{equation}\label{normal}  
\ttt \left(\nablas\cdot{\bfsig}_\Gamma(\nablas u_n)\right)\cdot\nablas=\bff\cdot \bfn
\end{equation}
and
\begin{equation}\label{tangential} 
 -{\bfsig}_\Gamma(\bfP_\Gamma\bfu)\cdot\nablas=\bfP_\Gamma\bff 
\end{equation}

\section{Plate Structures}

\subsection{Forces and Moments}

Consider first a subdomain polygonal subdomain $\omega \subset \Gamma$ 
of the plate $\Gamma$ with boundary $\partial \omega$ consisting of line 
segments $\gamma_i$. Using Greens formula on $\omega$ we obtain
\begin{align}\nonumber
&\ttt (\nablas \cdot ( \bfsig_\Gamma(\nablas u_n ) \cdot \nablas ), v_n )_\omega 
- (\bfsig_\Gamma (\bfu_t) \cdot \nablas ,\bfv_t )_\omega
\\ \label{eq-structures-a}
&\quad = -\ttt ( \bfsig_\Gamma(\nablas u_n )\cdot \nablas , \nablas v_n )_\omega 
+ (\bfsig_\Gamma(\bfu_t),\nablas \bfv_t )_\omega
\\ \nonumber
&\quad \quad + \ttt (\bfnu\cdot  (\bfsig_\Gamma(\nablas u_n ) \cdot \nablas) , v_n )_{\partial \omega} 
- ( \bfsig_\Gamma(\bfu_t) \cdot \bfnu , \bfv_t )_{\partial \omega}
\\ \label{eq-structures-b}
&\quad = \ttt (\bfsig_\Gamma(\nablas u_n ) , \bfeps_\Gamma(\nablas v_n) )_\omega 
+ (\bfsig_\Gamma(\bfu_t),\bfeps_\Gamma( \bfv_t ) )_\omega
\\ \nonumber
&\quad \quad + \ttt (\bfnu \cdot  (\bfsig_\Gamma(\nablas u_n )\cdot \nablas )\bfn - \bfsig_\Gamma(\bfu_t)\cdot \bfnu , \bfv )_{\partial \omega} 
\\ \nonumber
&\quad \quad - \ttt(\bfnu \cdot \bfsig_\Gamma(\nablas u_n )\cdot \bfnu, \nablas v_n)_{\partial \omega}
\end{align} 
where we used the identity $v_n = \bfv \cdot \bfn$ and moved the normal to the first slot in the bilinear form. Letting $\bftau$ be a unit tangent vector to $\partial \omega$, we may split the last term 
on the right hand side of (\ref{eq-structures-b}) in normal and tangent contributions as follows
\begin{align}\nonumber
&( \bfsig_\Gamma(\nablas u_n ) \cdot \bfnu , \nablas v_n)_{\partial \omega}
\\  \label{eq-structures-c}
&\quad =( \bfnu \cdot \bfsig_\Gamma(\nablas u_n ) \cdot \bfnu , \bfnu \cdot \nablas v_n)_{\partial \omega}
\\ \nonumber
&\quad \quad +
(\bftau \cdot \bfsig_\Gamma(\nablas u_n ) \cdot \bfnu , \bftau \cdot \nablas v_n)_{\partial \omega}
\end{align}
where the first term is the bending moment. For the second term on the right hand 
side (\ref{eq-structures-c}), integrating by parts along one of the line segments 
$\gamma_i$, with unit tangent and normal $\bftau_i = \bftau|_{\gamma_i}$ 
and $\bfn_i = \bfn |_{\gamma_i}$, we obtain 
\begin{align}\nonumber
&(\bftau_i \cdot \bfsig_\Gamma(\nablas u_n ) \cdot \bfnu_i , \bftau_i \cdot \nablas v_n)_{\gamma_i}
\\
&\quad = -( \bftau_i  \cdot \nablas (\bftau_i \cdot \bfsig_\Gamma(\nablas u_n ) \cdot \bfnu_i), v_n)_{\gamma_i}
\\ \nonumber
&\qquad \qquad 
+ (\bftau_i \cdot \bfsig_\Gamma(\nablas u_n ) \cdot \bfnu_i , v_n)_{\partial \gamma_i}
\\
&\quad = -( \bftau_i \cdot \nablas ( \bftau_i \cdot \bfsig_\Gamma(\nablas u_n ) \cdot \bfnu_i)\bfn, \bfv)_{\gamma_i}
\\ \nonumber
&\qquad \qquad 
+ (\bftau_i \cdot \bfsig_\Gamma(\nablas u_n ) \cdot \bfnu_i \bfn ,\bfv)_{\partial \gamma_i}
\end{align}
where $\partial \gamma_i$ consists of the two end points of the line segment $\gamma_i$. We 
introduce the following notation
\begin{align}
\label{eq:F}
\bfF &= \bfF_n + \bfF_t
\\
\label{eq:Fn}
\bfF_n &= \ttt  \bfnu \cdot  (\nablas \cdot \bfsig_\Gamma(\nablas u_n ))\bfn
\\ \nonumber 
&\quad - \ttt \bftau_i \cdot \nablas ( \bftau_i \cdot \bfsig_\Gamma(\nablas u_n ) \cdot \bfnu_i)\bfn
\\ \label{eq:Ft}
\bfF_t &= -  \bfsig_\Gamma(\bfu_t)\cdot \bfnu
\\ \label{eq:M}
M &=  \ttt \bfnu \cdot \bfsig_\Gamma(\nablas u_n ) \cdot \bfnu
\end{align}
for the normal and tangent components of the force and the moment at each 
of the line segments $\gamma$ on $\partial \omega$. Furthermore, we introduce 
the corner, or Kirchhoff, forces 
\begin{align}\label{eq:Fx}
\bfF_{\bfx,i} &= \bftau_i \cdot \bfsig_\Gamma(\nablas u_n ) \cdot \bfnu_i \bfn |_{\bfx}
\end{align}
at a corner $\bfx$ associated with a line segment $\gamma_i$, which has $\bfx$ as 
one of its endpoints and $\bftau_i$ is the unit tangent vector to $\gamma_i$ directed 
into $\bfx$. We then have the identity 
\begin{align}\nonumber
&\ttt (\nablas \cdot ( \bfsig_\Gamma(\nablas u_n ) \cdot \nablas ), v_n )_\omega 
- (\bfsig_\Gamma (\bfu_t) \cdot \nablas ,\bfv_t )_\omega
\\
&\quad 
= \ttt (\bfsig_\Gamma(\nablas u_n), \bfeps_\Gamma(\nablas v_n)_\omega 
+ (\bfsig(\bfu_t),\bfeps(\bfv_t))_\omega
\\  \nonumber
&\qquad +(\bfF,\bfv)_{\partial \omega}
- (M,\bfnu \cdot \nablas v_n )_{\partial \omega} 
+ \sum_{\bfx \in \mcX(\partial \omega)}  \sum_{i \in \mcI(x)} \bfF_{\bfx,i} 
\end{align}
where $\mcX(\partial \omega)$ is the set of corners on the polygonal boundary 
$\partial \omega$ and $\mcI(\bfx)$ is an enumeration of the two linesegments 
that has $\bfx$ as one of its endpoints.

\subsection{Jumps and Averages}

Consider a line segment $\gamma$ shared by two plates $\Gamma^+$ and $\Gamma^-$. 
We note that the force $\bfF^{\pm}$ is an $\IR^3$ valued 1-form  in $\bfnu^{\pm}$ and 
the moment $M^\pm$ is an $\IR$ valued 2-form in $\bfnu^\pm$.  More generally let 
$w^\pm = w^\pm(\bfnu^\pm,\dots \bfnu^\pm)$ be an $\IR^n$ valued $d$-linear form 
in $\bfnu^\pm$. Then we define the jump and average at $\gamma$ by
\begin{equation}\label{def-jump-average}
[w] = w^+ - (-1)^d w^-,\quad \langle w \rangle = \frac{1}{2} ( w^+ + (-1)^d w^- )
\end{equation} 
Note that when both plates $\Gamma^+$ and $\Gamma^-$ reside in the same plane 
$\bfnu^- = -\bfnu^+$ and we recover, using linearity and the simplified notation 
$w^\pm(\bfnu^\pm,\dots,\bfnu^\pm) = w(\bfnu^\pm)$, the standard jump 
\begin{align}
[w(\bfnu)] &= w^+(\bfnu^+) - (-1)^d w^-(\bfnu^-) 
\\
&= w^+(\bfnu^+) - (-1)^{2d} w^-(\bfnu^+)
\\
&= w^+(\bfnu^+) - w^-(\bfnu^+)
\end{align}
and similarly for the average. Finally, let $w^\pm_i$ be an $\IR^n$ valued $d_i$-linear form in 
$\bfnu^\pm$, then we note that $(w_1 \cdot w_2)^\pm = w_1^\pm \cdot w_2^\pm$ is an 
$\IR$ valued $(d_1 + d_2)$-linear form in $\bfnu^\pm$ and we have the identity
\begin{equation}\label{eq:jump-split}
[w_1 \cdot w_2 ] = [w_1] \cdot \langle w_2 \rangle + \langle w_1 \rangle \cdot  [w_2 ]
\end{equation} 
where for $n=1$ the scalar product is just usual multiplication of scalars. We may verify 
(\ref{eq:jump-split}) by 
\begin{align}
[w_1 \cdot w_2 ] &= w^+_1 \cdot w^+_2 - (-1)^{(d_1+d_2)} w^-_1 \cdot w^-_2
\\
&= w^+_1 \cdot w^+_2 - (-1)^{d_1} w^-_1  \cdot  (-1)^{d_2} w^-_2 
\\
&= w^+_1 \cdot w^+_2 - \widetilde{w}^-_1\cdot  \widetilde{w}^-_2
\\
&= (w^+_1 -  \widetilde{w}^-_1 ) \cdot \frac{w^+_2 +  \widetilde{w}^-_2}{2} 
\\ \nonumber
&\quad \quad + \frac{w^+_1 +  \widetilde{w}^-_1}{2} \cdot (w^+_2 - \widetilde{w}^-_2) 
\\
&= 
[w_1] \cdot \langle w_2 \rangle + \langle w_1 \rangle \cdot [ w_2 ]
\end{align}

%

\subsection{Interface Conditions}

Consider now a plate structure consisting of a finite number of plates such that at 
most two plates intersect in a common line segment. For simplicity we consider 
clamped boundary conditions on the boundary of the structure and focus our 
attention on the interface conditions at the intersections between the plates. 
For each line segment $\gamma$ where two plates $\Gamma^+$ and $\Gamma^-$ 
intersect we have the interface conditions 
\begin{align}
0&=[ \bfu ] \label{eq:interface-cont-u}
\\
0&=[\bfnu \cdot \nablas u_n  ] \label{eq:interface-cont-nnablau}
\\
0&=[\bfF] \label{eq:interface-F}
\\
0&=[M] \label{eq:interface-M}
\end{align}
corresponding to continuity of displacements, continuity of the rotation angle, 
equilibrium of forces, and equilibrium of moments. 

Furthermore, at each corner $\bfx$, not residing on the boundary of the structure, we 
require equilibrium of the Kirchhoff forces
\begin{equation}\label{eq:interface-Fx}
0 = \sum_{i \in \mcI(\bfx)} \bfF^+_{\bfx,i} +  \bfF^-_{\bfx,i} 
\end{equation}
where $\mcI(\bfx)$ is an enumeration of the line segments that meet in the corner $\bfx$ 
and $\bfF^\pm_{\bfx,i}$ is the Kirchhoff force emanating from plate $\Gamma_i^\pm$, 
the two plates that meet in line segment $i$. In other words, there are two contributions 
associated with each line segment, one for each of the two plates that share the line 
segment.

\section{Finite Element Formulation}

\subsection{The Mesh and Finite Element Space}

Let $\widehat{K}\subset \IR^2$ be a reference triangle and 
let $P_{2}(\widehat{K})$ be the space of polynomials of 
order less or equal to $2$ defined on $\widehat{K}$. Let 
$\Gamma$ be triangulated with quasi uniform triangulation 
$\mcK_{h}$ and mesh parameter $h\in (0,h_0]$ 
such that each triangle $K=F_{K}(\widehat{K})$ is planar
(a subparametric formulation). We let 
$\mcE_{h}$ denote the set of edges in the triangulation.

We here extend the discontinuous Galerkin method of Dedner et al. 
\cite{DeMaSt13} for the Laplace--Beltrami operator to the case of the plate.
We recall that $\Gamma$ is piecewise planar and thus
$\bfn$ is a piecewise constant exterior unit normal to $\Gamma$.

For the parametrization of $\Gamma$ we wish to define a map from 
a reference triangle $\widehat K$ defined in a local coordinate system 
$(\xi, \eta)$ to any given triangle $K$ on $\Gamma$. Thus the coordinates 
of the discrete surface are functions of the reference coordinates inside each element,
$\bfx_{\Gamma} = \bfx_{\Gamma}(\xi, \eta)$.
For any given parametrization, we can extend it to $\Omega_t$ by
defining
\begin{equation}\label{xdef}
\bfx(\xi,\eta,\zeta) := \bfx_{\Gamma}(\xi, \eta)+\zeta\,\bfn(\xi,\eta)
\end{equation}
where $-t/2 \leq \zeta \leq t/2$ and $\bfn$ is the normal to $\Gamma$.

We consider in particular a finite element parametrization of $\Gamma$ as
\begin{equation}\label{surfpara}
\bfx_{\Gamma}(\xi,\eta)= \sum_i\bfx_i\psi_i(\xi,\eta)
\end{equation}
where $\bfx_i$ are the physical location of the (geometry representing) nodes on the initial midsurface and $\psi_i(\xi,\eta)$ are affine finite element shape functions
on the reference element. (This parametrization is of course exact in the case of a piecewise planar $\Gamma$.)

For the approximation of the displacement, we use a constant extension,
\begin{equation}
\bfu \approx \bfu^h = \sum_i\bfu_i\varphi_i(\xi,\eta)
\end{equation}
where $\bfu_i$ are the nodal displacements, and $\varphi_i$ are piecewise quadratic shape functions.
We employ the usual finite element approximation of the physical derivatives of the chosen basis $\{\varphi_i\}$ on the surface, at $(\xi,\eta)$, in matrix representation, as
\begin{equation}
\left[\begin{array}{>{\displaystyle}c}
\frac{\partial \varphi_j}{\partial x}\\[2mm]
\frac{\partial \varphi_j}{\partial y}\\[2mm]
\frac{\partial \varphi_j}{\partial z}\end{array}\right] = \bfJ^{-1}(\xi,\eta,0) \left[\begin{array}{>{\displaystyle}c}
\frac{\partial \varphi_j}{\partial \xi}\\[2mm]
\frac{\partial \varphi_j}{\partial \eta}\\[2mm]
\frac{\partial \varphi_j}{\partial \zeta}\end{array}\right]_{\zeta=0} =: \bfJ^{-1}(\xi,\eta,0)\nabla_{\bfxi}\varphi_j\vert_{\zeta=0}
\end{equation}
where $\bfJ(\xi,\eta,\zeta) := \nabla_{\bfxi}\otimes\bfx$. This gives, at $\zeta=0$,
\begin{equation}\label{jacobiinv}
\left[\begin{array}{>{\displaystyle}c}
\frac{\partial \varphi_i}{\partial x}\\[2mm]
\frac{\partial \varphi_i}{\partial y}\\[2mm]
\frac{\partial \varphi_i}{\partial z}\end{array}\right] = \bfJ^{-1}(\xi,\eta,0) \left[\begin{array}{>{\displaystyle}c}
\frac{\partial \varphi_i}{\partial \xi}\\[2mm]
\frac{\partial \varphi_i}{\partial \eta}\\[2mm]
0\end{array}\right] 
\end{equation}
By (\ref{xdef}) we explicitly obtain
\begin{equation}
\left.\frac{\partial \bfx}{\partial \zeta}\right|_{\zeta=0}=\bfn 
\end{equation}
so
\begin{equation}
\bfJ(\xi,\eta,0) := \left[\begin{array}{>{\displaystyle}c>{\displaystyle}c>{\displaystyle}c}
\frac{\partial x}{\partial \xi} & \frac{\partial y}{\partial \xi} & \frac{\partial z}{\partial \xi}\\[3mm]
\frac{\partial x}{\partial \eta} & \frac{\partial y}{\partial \eta} & \frac{\partial z}{\partial \eta}\\[3mm]
n_x & n_y & n_z\end{array}\right] 
\end{equation}
We can now introduce finite element spaces constructed from the basis previously discussed by defining
\begin{align}\nonumber
W^h := {}& \{ v: {v\vert_T \circ F_K\in P^2(\widehat{K}),\; \forall K\in\mathcal{K}_h};\\ {}& v\; \in C^0(\Gamma),\;
  v=0\,\text{on $\partial\Gamma_\text{D}$}\} \label{spacevA}
\end{align}
We also need the set of interior edges defined by
\begin{align}
  \mcE_h^\text{I}&:=\{E = K^+ \cap K^-:  K^+, K^- \in \mcK_h \}
\end{align}
and the set of boundary edges on the Dirichlet part of the boundary 
\begin{align} 
  \mcE_h^\text{D} &:= \{ E =  K \cap \partial\Gamma_\text{D}: \; K \in \mcK_{h} \}
\end{align}
To each interior edge $E$ we associate the conormals $\bfnu^{\pm}_E$ 
given by the unique unit vector which is tangent  to the surface element 
$K^{\pm}$, perpendicular to $E$ and points outwards with respect to 
$K^{\pm}$. 
Note that the conormals $\bfnu_E^{\pm}$ may lie in different planes 
at junctions between different plates. The jump and average of 
multilinear forms for edges $E \in \mcE_h^I$ are defined by 
(\ref{def-jump-average}). For edges $E \in \mcE_h^D$ it is convenient 
to use the notation 
\begin{equation}
\langle w \rangle = [w] = w
\end{equation}


\subsection{The Method}

Our finite element method takes the form: find $\bfU \in \bfV_h := [W_h]^3$ 
such that
\begin{equation}\label{eq:fem}
A_h(\bfU,\bfv) = l_h(\bfv) \quad \forall \bfv \in \bfV_h 
\end{equation}
Here the bilinear form $A_h(\cdot, \cdot)$ is defined by
\begin{align}
A_h(\bfv,\bfw) &:= a_h^\text{P}(\nabla_{\Gamma} v_n,\nabla_{\Gamma} w_n )
+ a_h(\bfv_t,\bfw_t) \label{eq:BilinearMembrane}
\end{align}
with  $\bfv = v_n \bfn + \bfv_t$  and
\begin{equation}
a_h(\bfv_t,\bfw_t) := \sum_{K\in\mcK_h}( \bm{\sigma}_{\Gamma}(\bm{v}_t),\bm{\varepsilon}_{\Gamma}(\bm{w}_t))_{K} 
\end{equation}
where $(\cdot,\cdot)_\omega$ denotes the $L_2(\omega)$ scalar product, and
\begin{align}
a_h^\text{P}(\bfv,\bfw) := {}&  \ttt a_h(\bfv_t,\bfw_t) 
\\ \nonumber
&{} -\sum_{E \in \mcE_h^\text{I}\cup\mcE_h^\text{D}}(\mean{M(\bfv)},\jump{\bfnu_E\cdot\bfw})_E 
\\ \nonumber
&{} -\sum_{E \in \mcE_h^\text{I}\cup\mcE_h^\text{D}}(\mean{M(\bfw)},\jump{\bfnu_E\cdot\bfv})_E 
\\ \nonumber
&{}  +\frac{\beta \ttt}{h}\sum_{E \in \mcE_h^\text{I}\cup\mcE_h^\text{D}}(\jump{\bfnu_E\cdot\bfv},\jump{\bfnu_E\cdot\bfw})_E
\end{align}
Here $\beta=\beta_0(2\mu+2\lambda)$ where $\beta_0$ is an $O(1)$ constant, 
cf. \cite{HaLa11}, and we also recall that the factor $\ttt$ is included in the definition 
(\ref{eq:M}) of the moment $M$.  The right hand side is given by 
\begin{equation}
l_h(\bfv) := ( \bff, \bfv )_{\Gamma}
\end{equation}

This is a c/dG method closely related to the one studied in \cite{HaLa11}, with the difference 
of being formulated in an arbitrary orientation in $\mathbb{R}^3$, including membrane 
deformations, and extended to structures of plates. 

We note that: 
\begin{itemize}
\item The continuity of displacement (\ref{eq:interface-cont-u}) is strongly 
enforced since $\bfV_h$ consists of continuous functions. 
\item The continuity of the rotation angle (\ref{eq:interface-cont-nnablau}) is weakly 
enforced by the discontinuous Galerkin method. 
\item The force equilibrium conditions (\ref{eq:interface-F}) and (\ref{eq:interface-Fx}) 
are weakly enforced but does not give rise to any additional terms in the formulation since 
$\bfV_h$ consists of continuous functions. 
\item The moment equilibrium condition (\ref{eq:interface-M}) is weakly enforced 
by the discontinuous Galerkin method. 
\end{itemize}
More precisely, consider an edge $E \in  \mcE^I_h$ shared by two elements $K^+$ and $K^-$. 
Multiplying the exact equation by a test function $\bfv \in \bfV_h$ and using Green's 
formula element wise generates the following contribution at the edge $E$,
\begin{align}
&(\bfF^+,\bfv^+)_\gamma + (\bfF^-,\bfv^-)_E 
\\
&\quad - 
(M^+, \bfnu^+_E \cdot v_n^+)_E + (M^+, \bfnu^+ \cdot v_n^+)_E
\end{align}
where $\bfF^\pm = \bfF^\pm(\bfu)$ and  $M^\pm = M^\pm(\bfu)$. For the first term 
we have using the continuity of $\bfv$ and (\ref{eq:F}), 
\begin{equation}
(\bfF^+,\bfv^+)_E + (\bfF^-,\bfv^-)_E 
= 
([\bfF],\bfv)_E
=
0
\end{equation}
For the second term we note that the integrand may be written 
\begin{align}
M^+ \bfnu^+_E \cdot v_n^+ + M^+ \bfnu^+ \cdot v_n^+ 
= [ M   \bfnu \cdot v_n ]
\end{align}
where we used the fact that $M^\pm$ is $2$-linear in $\bfnu^\pm$, see (\ref{eq:M}), 
and $\bfnu^\pm \cdot \nablas v^\pm_n$ is 1-linear in $\bfnu$, and thus 
$M^\pm \bfnu^\pm \cdot \nablas v^\pm_n$ is 3-linear in $\bfnu^\pm$, together with the definition 
(\ref{def-jump-average}) of the jump to write the sum as a jump. Next using 
(\ref{eq:jump-split}) we get
\begin{align}
[M \bfnu \cdot \nablas v_n ] 
&=
[M] \langle \bfnu \cdot \nablas v_n \rangle 
+ 
\langle M \rangle [ \bfnu \cdot \nablas v_n ]
\\
&= \langle M \rangle [ \bfnu \cdot \nablas v_n ]
\end{align}
since $[M] = 0$ according to (\ref{eq:interface-M}). Thus the second term takes the form
\begin{align}
&- (M^+(\bfu), \bfnu^+_E \cdot \nablas v_n^+)_E 
- (M^+(\bfu), \bfnu^+ \cdot \nablas v_n^+)_E
\\
&\quad = - ( \langle M(\bfu) \rangle ,  [ \bfnu \cdot \nablas v_n ] )_E
\\
&\quad=  -( \langle M(\bfu) \rangle, [ \bfnu \cdot \nablas v_n ] )_E 
\\ \nonumber
&\quad \quad \quad - ( \langle M(\bfv) \rangle , [ \bfnu \cdot \nablas u_n ])_E
\end{align}
where at last we symmetrized using the fact that the added term is zero by 
(\ref{eq:interface-cont-nnablau}) and we included the dependency $M = M(\bfu)$ for clarity. 
We finally note that we have the following identities
\begin{equation}
\langle M \rangle = \frac{1}{2} (M^+ + M^-)
\end{equation}
and 
\begin{equation}
[ \bfnu \cdot \nablas v_n ] =  \bfnu^+ \cdot \nablas v^+_n +   \bfnu^- \cdot \nablas v^-_n 
\end{equation}

\begin{rem}
We note that the method for a plate structure has the same form as for a single plate since 
we use the proper definitions of jumps and averages encoded by the conormal.
\end{rem}

\begin{rem}
We note that with this formulation, we have Galerkin orthogonality
\begin{equation}
A_h(\bfu-\bfU,\bfv) = 0 \quad \forall \bfv \in \bfV_h
\end{equation}
which enables us to prove an a priori error estimate of optimal order provided the solution is regular 
enough using the same techniques as in \cite{HaLa02}.
\end{rem}
\begin{rem}
For shell modelling, the plate approach can still be used by viewing the shell as an assembly of facet elements. Then we have an elementwise planar approximation $\Gamma_h$ of $\Gamma$ and we use
elementwise projections $\bfP_h =\bfI-\bfn_h\otimes\bfn_h$, where $\bfn_h$ is the elementwise constant approximation of $\bfn$. The differential operators are then defined on the
discrete surface, e.g,
$\nablash v := \bfP_h\nabla v$, etc., and replacing the exact differential operators and exact surface by their discrete approximations in (\ref{eq:fem}) we obtain a simple shell model.
\end{rem}
\section{Numerical Examples}

\begin{figure}[ht]
\centering
\includegraphics[width=10cm]{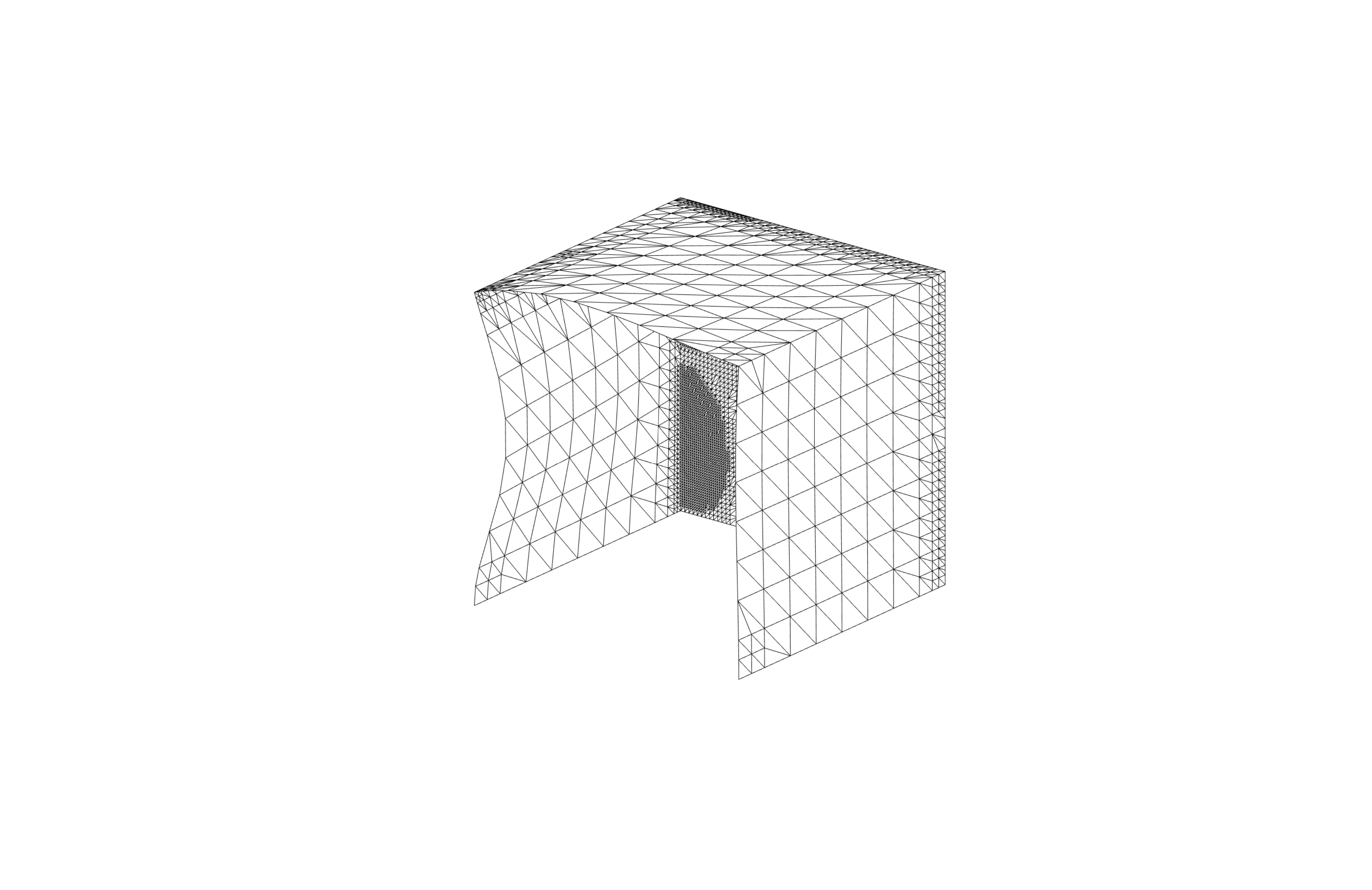}
\includegraphics[width=10cm]{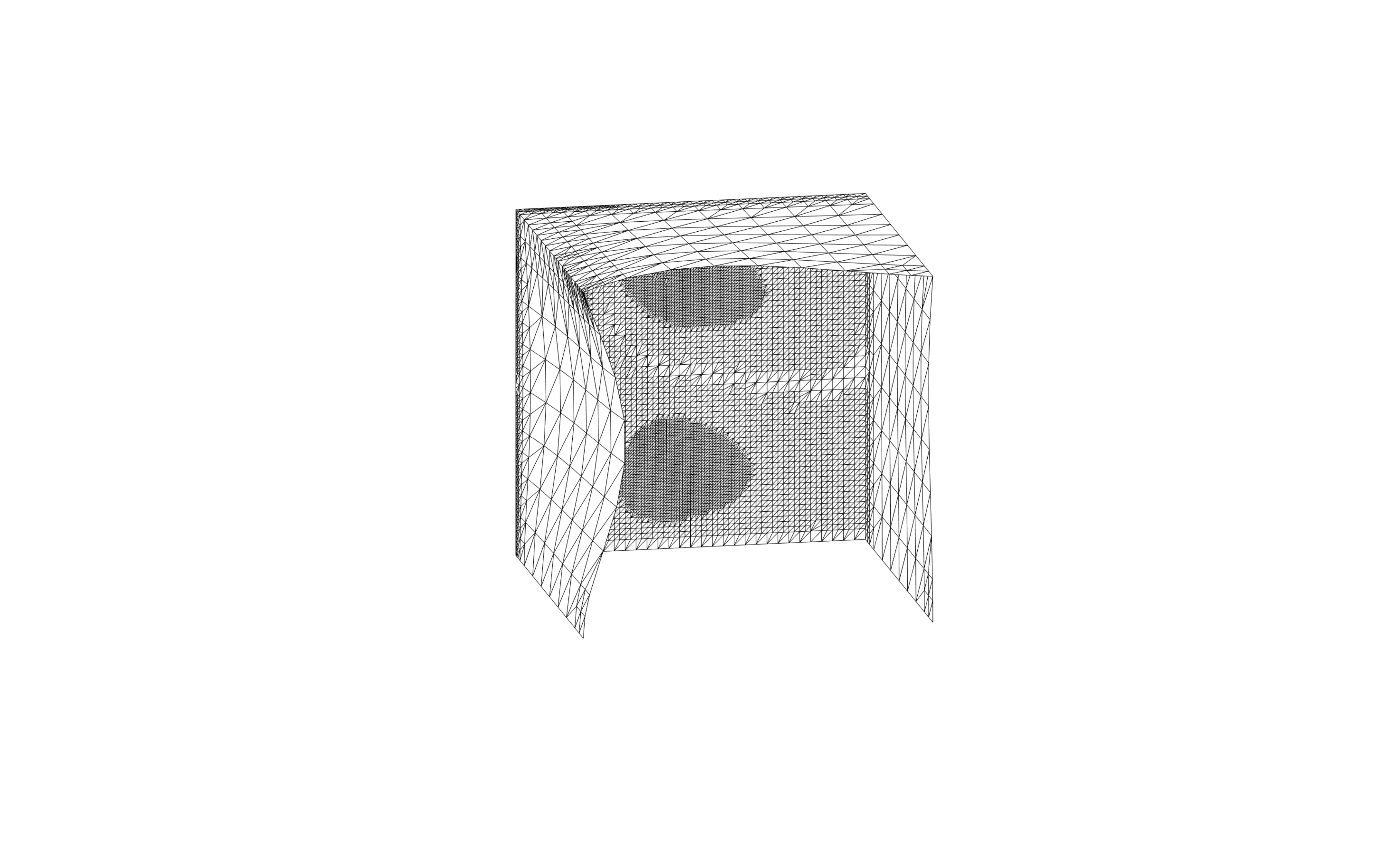}
\includegraphics[width=8cm]{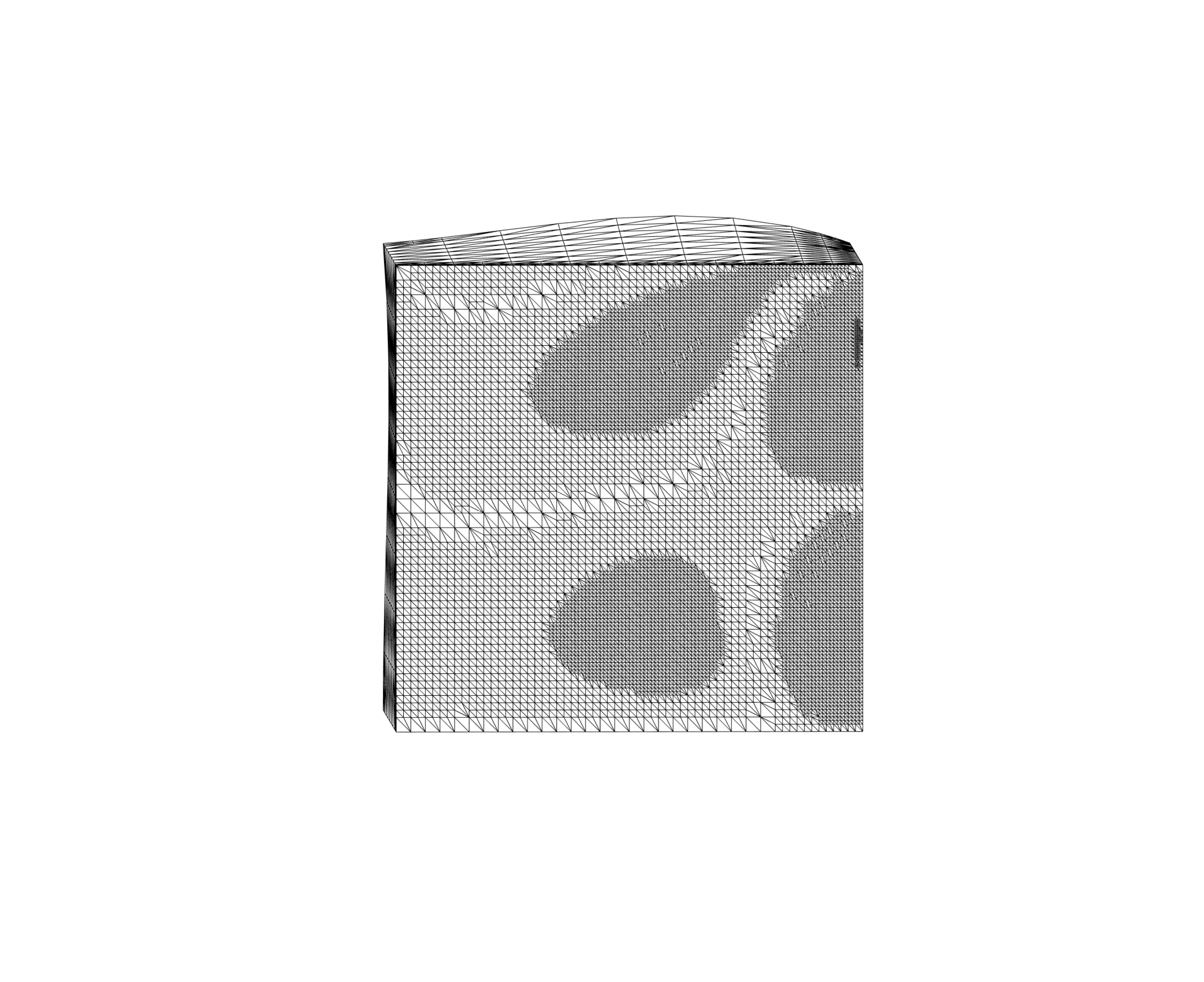}
\caption{Different views of the deformed box with $t=10^{-3}$.}\label{tminus3}
\end{figure}

\begin{figure}[ht]
\centering
\includegraphics[width=10cm]{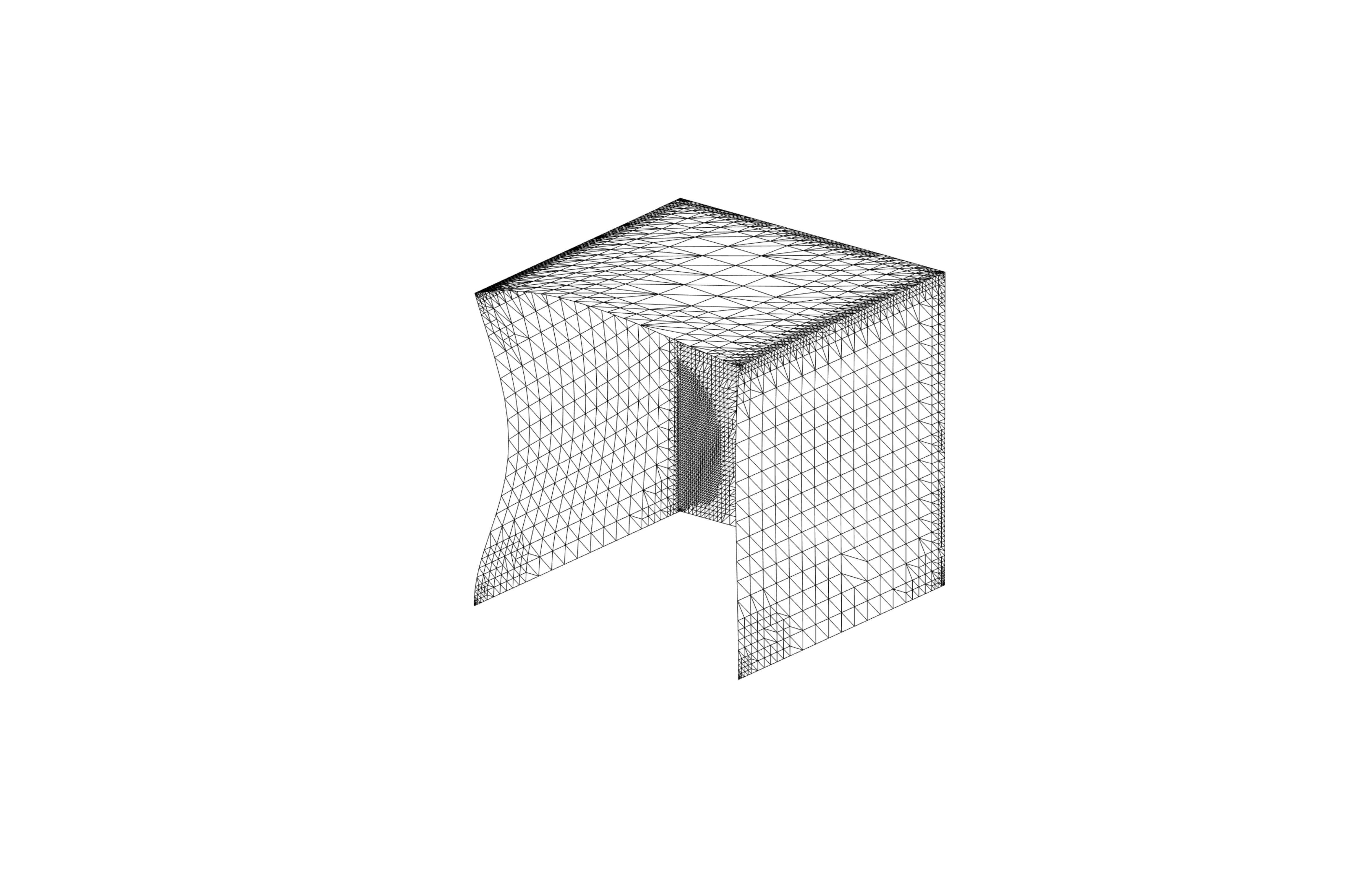}
\includegraphics[width=10cm]{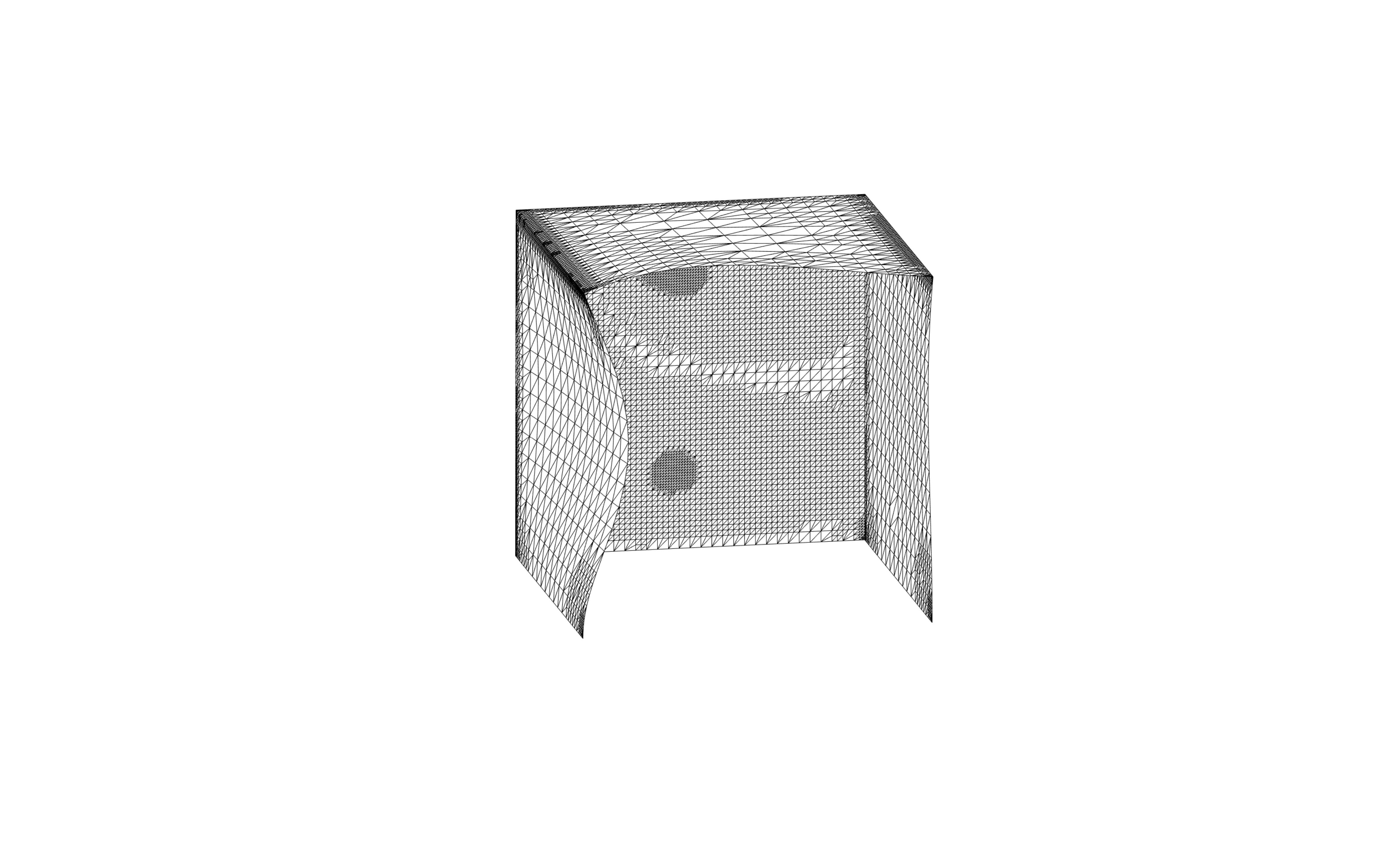}
\includegraphics[width=8cm]{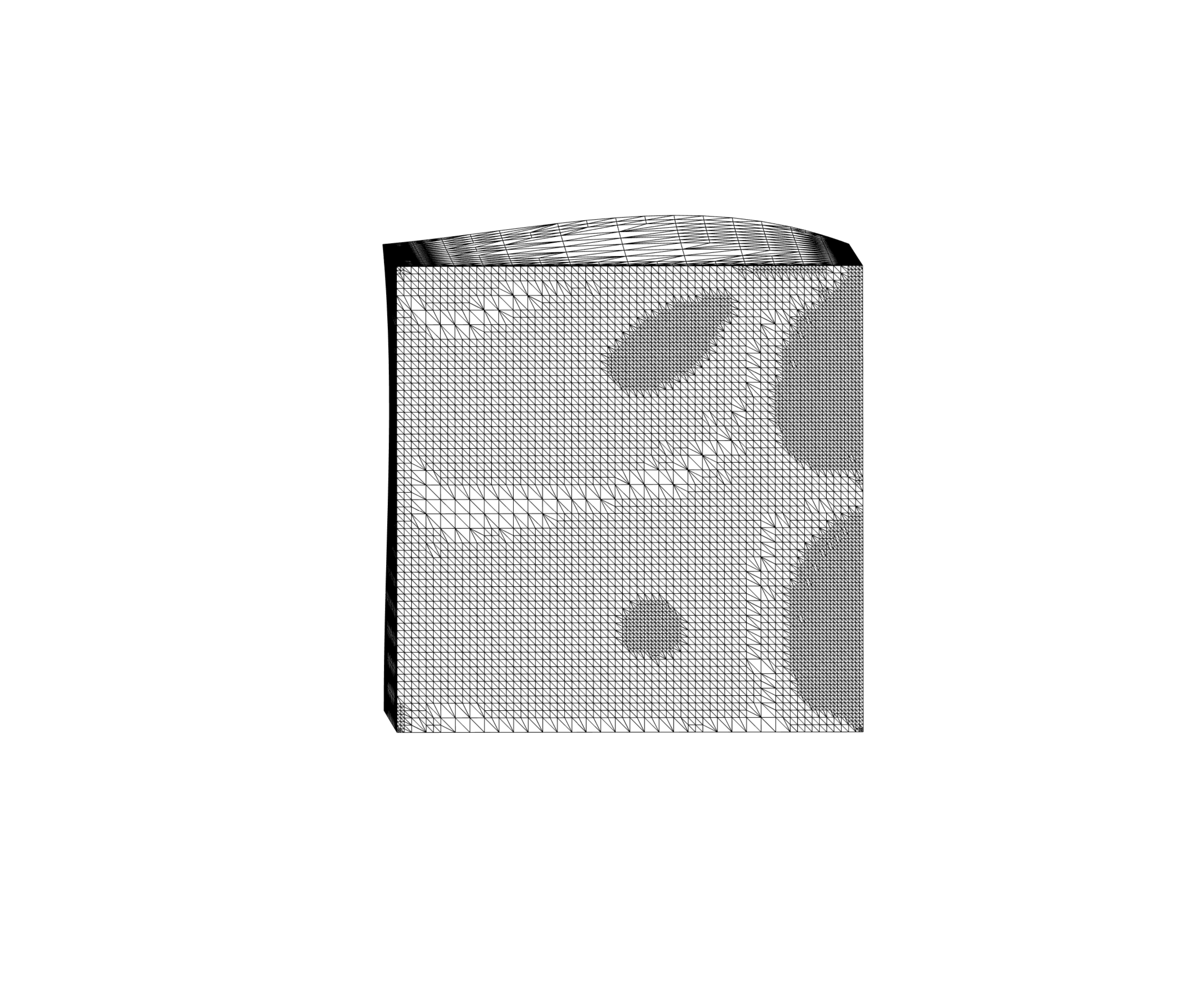}
\caption{Different views of the deformed box with $t=10^{-2}$.}\label{tminus2}
\end{figure}

\begin{figure}[ht]
\centering
\includegraphics[width=10cm]{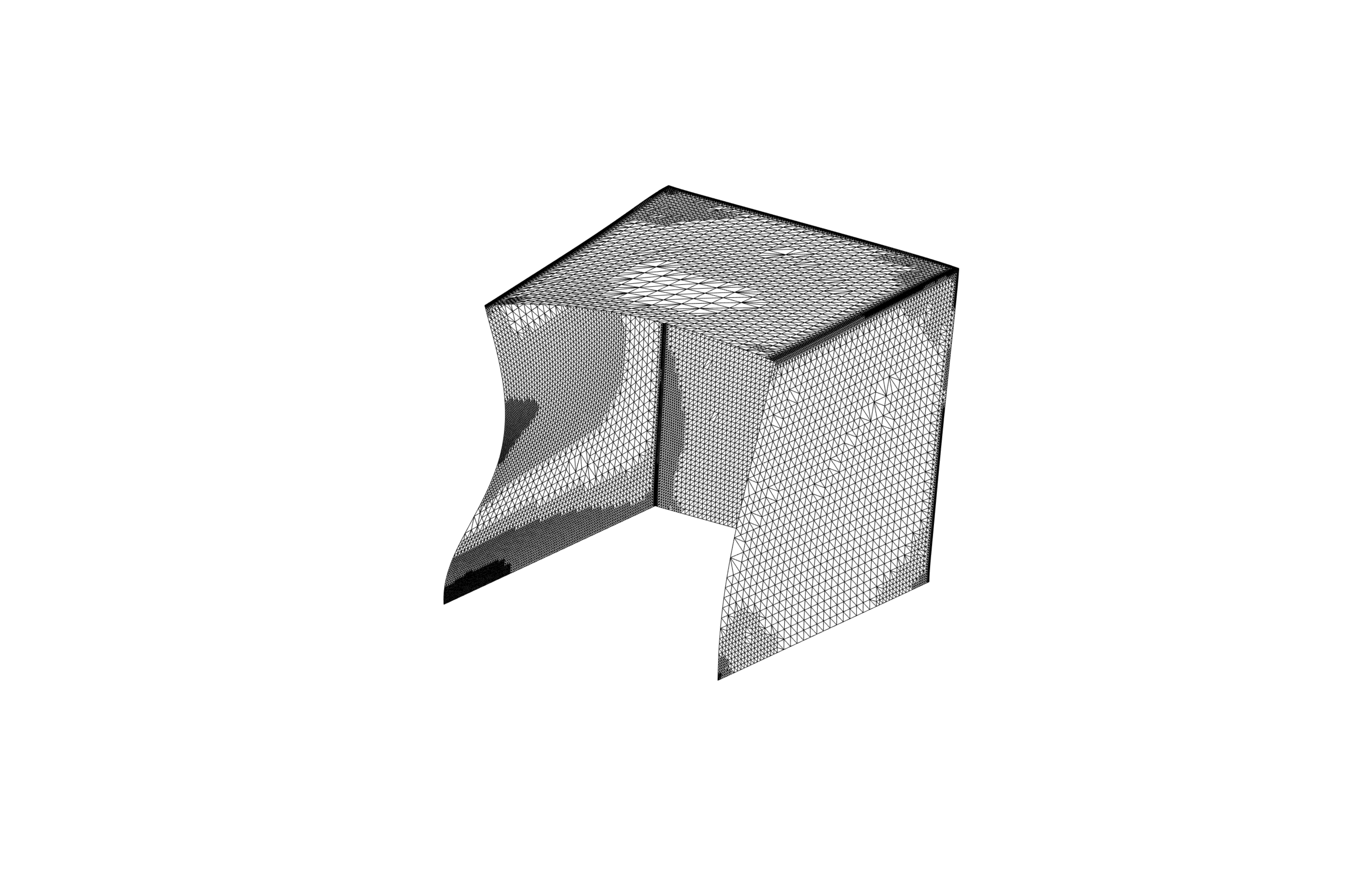}
\includegraphics[width=10cm]{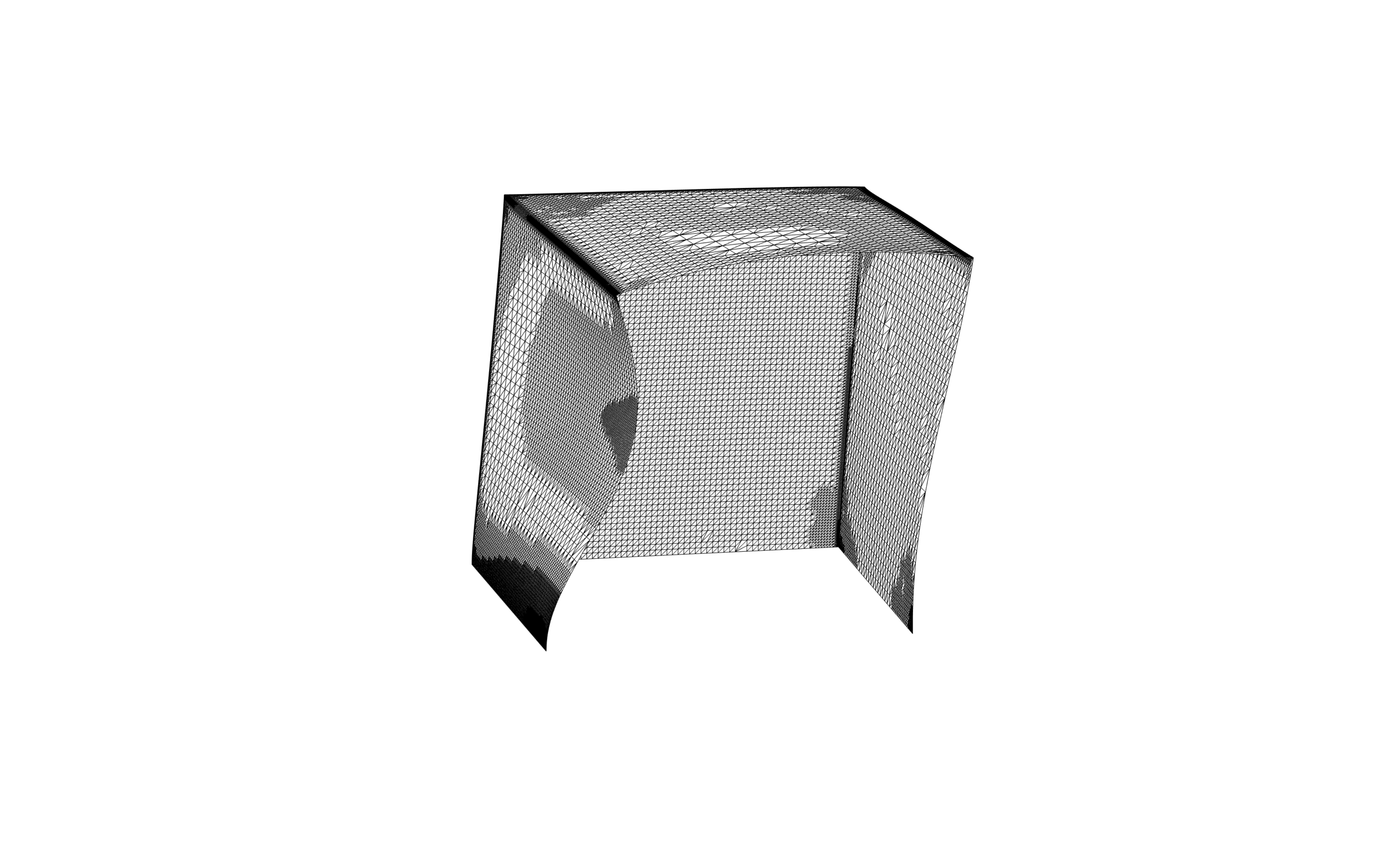}
\includegraphics[width=8cm]{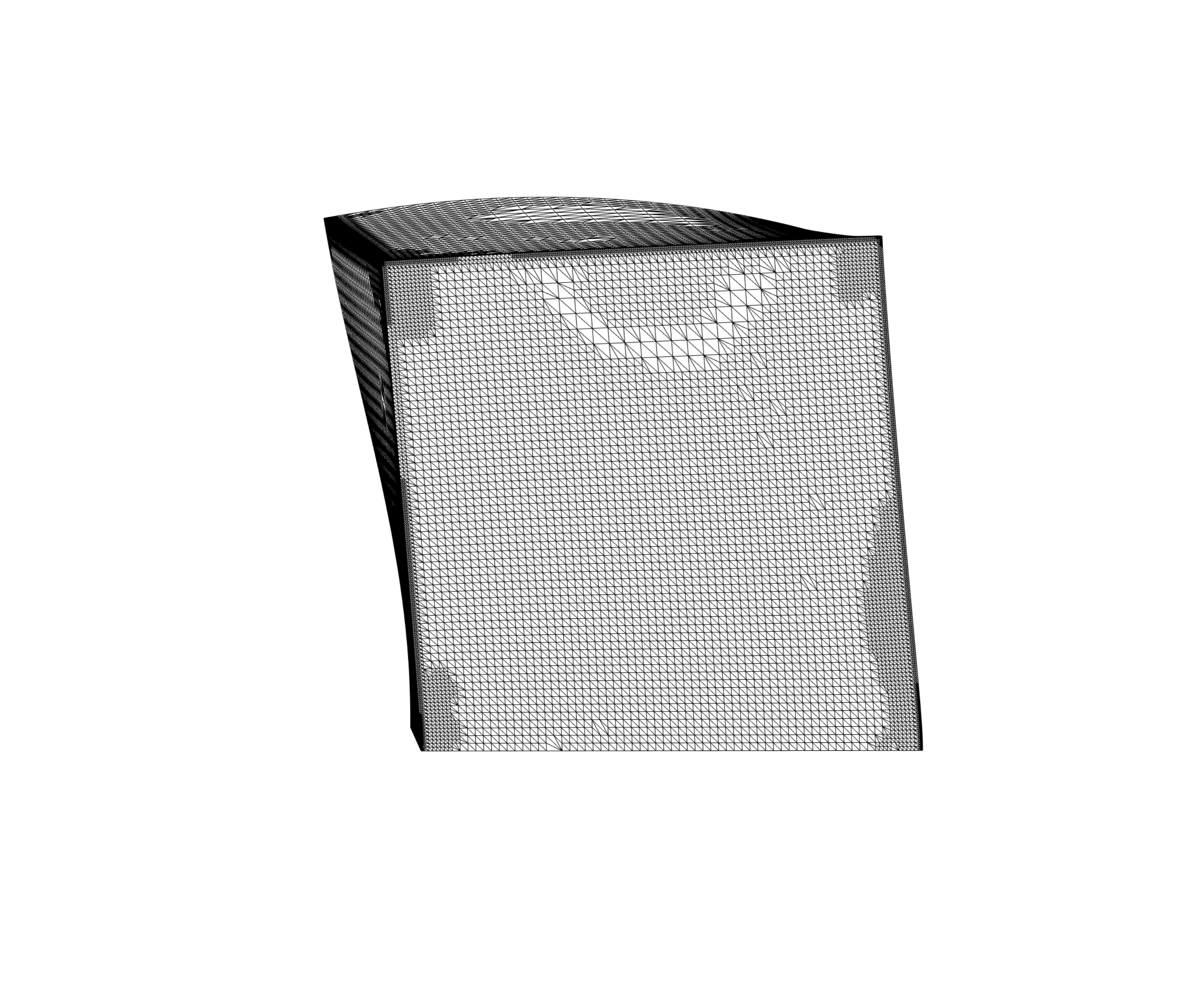}
\caption{Different views of the deformed box with $t=10^{-1}$.}\label{tminus1}
\end{figure}

We consider the surface of the box $[0,1]\times [0,1] \times [0,1]$, fixed to the floor and with one wall missing. The material data are: Poisson's ratio $\nu = 0.5$ and Young's modulus $E = 10^9$. The stabilization parameter was set to $\beta_0 = 10$.
An \emph{ad hoc}\/ residual--based adaptive scheme was used to generate locally refined meshes. The load was given as
\[
\bff=t^2\left[\begin{array}{c}
4\times 10^7\\
0\\
0\end{array}\right]
\]
at $x=0$, $\bff =\boldsymbol 0$ elsewhere. The point of the scaling with thickness is that 
after division by $t^2$ the membrane stiffness will scale with $t^{-2}$ so that the limit of $t\rightarrow 0$ corresponds to the inextensible plate solution. With increasing $t$ the membrane effect will become more and more visible. The numerical results using three different thicknesses, $t = 10^{-k}$, $k=3,2,1$, are given in Figs. \ref{tminus3}--\ref{tminus1}. Note the marked membrane deformations at $k=1$.

\section{Concluding Remarks}
In this paper we have introduced a c/dG method for arbitrarily oriented plate structures.
Our method  is expressed directly in the spatial coordinates, unlike traditional schemes that typically are based on coordinate transformations from planar elements.
This leads to a remarkably simple and easy to implement discrete scheme. The c/dG approach also allows for avoiding the use of $C^1$--continuity, otherwise required by the plate model, by allowing for discontinuous rotations between elements, and the same function space can then be used to model both plate and membrane deformations. We also introduced the proper conormals, mean values, and jumps necessary for handling the discontinuities on the element borders.

\begin{acknowledgements}
This research was supported in part by the Swe\-dish Foundation for Strategic Research Grant
No.\ AM13-0029, the Swedish Research Council Grants Nos.\ 2011-4992,
2013-4708, and Swedish strategic research programme eSSENCE.
\end{acknowledgements}

\newpage


\end{document}